\documentclass[a4paper,11pt]{article}

\usepackage{natbib}        
\usepackage[dvips]{graphicx} 
\usepackage[latin1]{inputenc}
\usepackage{amssymb}
\usepackage{a4wide}         
\usepackage{amsmath}         
\usepackage{amsfonts}
\usepackage[usenames,dvipsnames]{xcolor}
\usepackage{curves}

\begin{document}

\title{\LARGE \bf About string stability of a vehicle chain with unidirectional controller} 


\author{Arash Farnam\thanks{ID Lab, Department of Electronics and Information Systems (ELIS), Faculty of Engineering and Architecture, Ghent University; Technologiepark Zwijnaarde 914, 9052 Zwijnaarde(Ghent), Belgium (e-mail: arash.farnam@ugent.be).}~~and Alain Sarlette\thanks{ELIS, Ghent University, Belgium; and QUANTIC lab, INRIA Paris, France (e-mail: alain.sarlette@inria.fr)}} 

\maketitle

\begin{abstract}                
This paper deals with the problem of string stability in a chain of acceleration-controlled vehicles. It is known that string stability cannot be achieved, with any linear controller, when the vehicles' control inputs are based on relative distances to a fixed number of predecessors. We extend the set of impossible settings by including elements like dynamic sensor parts and local inter-vehicular communication, as in cooperative adaptive cruise control. It is also known that a weaker form of string stability is achievable by adding absolute velocity measurements (e.g.~``time-headway'' policy). We show that a stronger property can also be achieved, provided steady-state control gain is infinite e.g.~by using integral control.
\end{abstract}



\section{Introduction}\label{sec:intro}

Grouping vehicles into tight platoons is a method for e.g.~increasing the capacity of roads by automated highway systems (\cite{1}). The distances between vehicles is decreased by ensuring safety thanks to automatic controllers, enabling many vehicles to accelerate or brake simultaneously and eliminating the distance needed for human reaction. The most fundamental platoon is the vehicle chain, where all vehicles are aligned after each other. During the recent years numerous works have considered different control strategies to stabilize each vehicle at a desired distance from its predecessor and follower (\cite{1,2,3,4,Dirk}).

When disturbance inputs affect the vehicle chain, its a priori cooperative coupling can lead to new types of instabilities. In particular, \emph{string instability} is a situation where the spacing error between consecutive vehicles grows unbounded as the number of vehicles increases to infinity. Since its definition in \cite{6,8}, string (in)stability has spurred a lot of discussion and research. Indeed, it is well known since \cite{8} that string stability cannot be achieved in a homogeneous string of interconnected second-order integrators (e.g.~acceleration-controlled vehicles), with \emph{any} controller that is linear and whose local control actions are determined from the relative distance to a few directly preceding vehicles: some perturbations will unavoidably grow unbounded along the chain. This has attracted attention as a prototypical, unavoidable shortcoming of linear systems (\cite{3,4}). When the controller reacts to the just preceding vehicle only, this follows essentially from the Bode integral for the transfer function from vehicle $i-1$ to $i$, which takes the form of a complementary sensitivity function.

To solve this issue, a notable alternative setting is to allow the use of \emph{absolute} velocity in the controller, see \cite{Dirk}, \cite{15}, \cite{16}. It has been shown that the effect of a disturbance from the leading vehicle can then be kept in check for arbitrarily long chains. Another natural feature is to allow local communication, as in Cooperative Adaptive Cruise Control (CACC, see \cite{7}, \cite{17}, \cite{18}, \cite{19}, \cite{20}). Remarkably, the literature considering communication also always assumes the use of absolute velocity in the controller, and it thus remains unclear what can be done with communication but without using absolute velocity.

The aim of the present paper is to add essentially two types of precision to this picture. First, we consider countering the effect of disturbances acting possibly \emph{on all vehicles} in the $(L_2,l_2)$ sense. For this, we give a positive result using PID control and absolute velocity. Second, we clarify how alternative settings behave in absence of absolute velocity or integral action in the controller. For this, we provide a series of negative results; in particular the use of absolute velocity appears necessary even in presence of CACC-type communication. A more detailed statement of state of the art and of our contributions is given in Section \ref{sec:3}, after clarifying the setting.

While the present work focuses on unidirectional vehicle chains -- i.e.~vehicles react to their predecessors only -- a whole line of work has been developed for bidirectional chains as well, where vehicles react to predecessors and followers. Impossibilities to satisfy string stability have been obtained for both symmetric (\cite{9,11}) and asymmetric (\cite{23}) couplings. Conversely, \cite{12,13} have identified the possibility to avoid unbounded growth of a disturbance that acts \emph{on the leader only}, without resorting to absolute velocity; this result thus fundamentally differs from the impossibility implied by the Bode integral for unidirectional chains. Our viewpoint on the bidirectional setting takes a completely different approach and will be given elsewhere.


\section{Problem Setting}\label{sec:setting}


\subsection{Open-loop model and control objective}

Consider a chain of $N+1$ subsystems (e.g.~vehicles) with respective configuration (e.g.~deviation from nominal position) $x_i \in \mathbb{R}$, for $i=0,1,2,...,N$. They move according to the second-order integrator dynamics, expressed in Laplace domain:
\begin{equation}\label{eq:2order}
s^2 x_i(s) = u_i(s) + d_i(s)
\end{equation}
where $u_i$ is the control input, and $d_i$ is a disturbance input signal. These disturbance signals can also be an indirect way to model nonzero initial conditions, which are the focus in other papers. The pure second-order integrator is a standard idealization, valid in good approximation for space vehicles, vacuum transit, after other dynamics has been compensated by local feedback, or when it can be included (e.g.~by linearization) in the form of $u_i(s)$ that we further specify below.

The goal of stability, on a system like \eqref{eq:2order} with fixed $N$, is to ensure that arbitrary input signals $d_i$ are not amplified unboundedly in the state or output. \emph{String} stability further checks what happens when $N$ becomes infinite; thus if the system is stable for each $N$ but with stability bounds depending on $N$ in a bad way, then string stability may fail. The standard definition of string stability (\cite{6}, \cite{11}) considers how the relative distances between consecutive vehicles are affected. In the ``time-headway'' extension \cite{15}, the target distance between vehicles is made proportional to their absolute velocity, such that the configuration error (``output signal'' of interest) is
\begin{equation}\label{eq:errmeas}
e_i(s)= x_{i-1}(s) - x_{i}(s) - hs x_i(s) \, ,
\end{equation}
for some $h \geq 0$. For a vector $v$ of signals $v_i(t)$, $i=1,2,...,N$, the $L_2$ norm of $v_i$ is defined as
$\Vert v_i(.) \Vert_2 = \sqrt{\int_{-\infty}^{+\infty} \, \left(v_i(t)\right)^2 \, dt}$
and the $(L_2,l_2)$ norm of $v$ is defined as
$\Vert v(.) \Vert_{2} =  \sqrt{\sum_{i=1}^N\int_{-\infty}^{+\infty} \, \left(v_i(t)\right)^2 \, dt} \; ,$
which by Parseval both have equivalent expressions in frequency domain.

\noindent \textbf{Definition 1 [$(L_2,l_2)$ string stability]:} \emph{The chain \eqref{eq:2order}, controlled with feedback signals $u_i$ to be designed, is called $(L_2,l_2)$ string stable if there exists a constant value $c_1$ such that 
$$\Vert e(.) \Vert_2 \leq c_1 \; \Vert d(.) \Vert_2 $$
for all bounded signals $d$, and all chain lengths $N$.}

The difficulty is to make the bound uniform in $N$. This definition, which can be found in \cite{2,9,11} among others, requires that the disturbance $d_k$ acting at a specific subsystem $k$, should have a \emph{strictly decreasing} influence along the vehicle chain; indeed, if each disturbance was just transported like a wave along the chain, then $e_N$ would effectively be subject to their sum, and thus unboundedly increasing with $N$. However, pure transport would satisfy a somewhat weaker requirement, considered in e.g.~\cite{10,12,13,15,16,23}. We would call it $L_2$ string stability.

\noindent \textbf{Definition 2 [$L_2$ string stability]:} \emph{The chain \eqref{eq:2order}, controlled with feedback signals $u_i$ to be designed, is called $L_2$ string stable if there exists a constant value $c_1$ such that 
$$\Vert e_k(.) \Vert_2 \leq c \; \Vert d_i(.) \Vert_2  \quad \text{for all } k \; ,$$
for any situation with bounded disturbance signal $d_i$ at some $i$ and zero disturbances on all $k\neq i$, and for all $N$.}

This definition essentially requires that the transfer function from any $d_i$ to any $e_k$ is $H$-infinity bounded \emph{independently of $N$}. It is necessary for $(L_2,l_2)$ string stability, but not sufficient.


\subsection{Constraints on controller design}

We will consider different constraints on the feedback controller design, all revolving around the unavailability of absolute or long-range position measurements. The most general controller that we consider, and motivated in the next paragraph, writes:
\begin{eqnarray}\label{gc0}
u_i(s) & = & K(s) e'_{i}(s) + H(s) r_i(s)\\
\nonumber r_{i}(s) & = & W(s)\, v_{i-1}(s) \;\; \\
\nonumber v_{i}(s) & = & F(s) e'_{i}(s) + G(s) r_i(s) \;\; \\
\nonumber e'_i & = & M^{(r)}(s) x_{i-1}(s) - M^{(f)}(s) x_i(s) - h s x_i(s) \; ,
\end{eqnarray}
for $i=1,2,...,N$. Here $K(s)$ represents the linear controller's transfer function depending on measured configuration error $e'_i$, and satisfying $K(0) \neq 0$ (no poles cancellation).  Furthermore, unlike in \cite{22}, we do not allow $K(s)$ to grow unboundedly with $N$ and we assume all transfer functions independent of $N$.
The $e'_i$ can differ from the configuration error $e_i$ through the transfer functions $M^{(r)}(s), M^{(f)}(s)$. They express that inter-vehicle distances are typically measured by sensors which could \emph{break the symmetry of perfect relative distances}. Indeed, consider distance sensors with parts mounted on each vehicle, and these mounts can have some dynamics. These dynamics are themselves sensitive to the relative position of a sensor part $p_i$ with respect to the vehicle on which it is mounted $x_i$. Thus, we should have $\; s^2 p_i = K^{(.)}(s)\, (x_i-p_i)\;$,
which yields:
\begin{eqnarray}
\label{eq:sens2}
M^{(r)} &=& \frac{K^{(r)}}{s^2 +K^{(r)}} \; ,\quad M^{(f)} = \frac{K^{(f)}}{s^2 +K^{(f)}} \; ,
\end{eqnarray}
with $K^{(r)}(0)\neq 0$ and $K^{(f)}(0) \neq 0$, respectively for the sensor parts mounted on the rear and front ends of the vehicles. In particular, at the limit of infinitely stiff mounts $K^{(r)},K^{(f)} \rightarrow \infty$, we get $$M^{(r)}(s) = M^{(f)}(s) = \text{Identity}$$ and the controller $u_i(s)$ just depends on the configuration error $e_i(s)$; if furthermore $h=0$ it reduces to relative position measurement $e_i(s)= x_{i-1}(s) - x_{i}(s)$. In addition, $u_i$ relies on a communication signal $r_i(s)$ received from the preceding vehicle. Imperfections in the communication channel are taken into account by a simplified linear model, with the bounded transfer function $\mid W(j\omega)\mid$ not exactly known. The signal $v_i(s)$ sent into this channel by vehicle $i$ is computed as a linear function, on the same basis as the control signal, with $H(s), F(s)$ and $G(s)$ controller transfer functions to be designed. We make the following assumptions on the communication channel.

\noindent \textbf{Assumption 1 [communication channel]:}\vspace{-3mm}
\begin{itemize}
\item (stability) The poles of $F,G,H,W$ must all have negative real parts.
\item (keeping unmodeled noise in check) $G(j\omega)$ and $F(j\omega)$ are bounded for all $\omega$, while $H$ may be of a similar form as $K$.
\item (robustness) The controller cannot rely on perfect cancellation effects by matching of $W$ with $F,G,H$.
\end{itemize}

Let us briefly motivate the various features of \eqref{gc0}. 
\newline - The most basic setting, with $h=0$, $M^{(r)}(s) = M^{(f)}(s) = \text{Identity}$ and $W=0$ (thus no communication $H=F=G=0$), was the initial focus of \cite{6} and is known to be string unstable with \emph{any linear controller}. To check how alternative settings could help, we stay in the linear realm. 
\newline - Just introducing the possibility of $h>0$, was shown to enable satisfying Definition 2 (\cite{15}). However, the time-headway policy makes the effective inter-vehicle distance velocity-dependent, where absolute velocity of the chain would be defined via some other controller, to be carefully interfaced. It also remains to be seen exactly how (accurately) the absolute velocity would be measured in practice towards implementing \eqref{gc0} with $h>0$. Therefore, it is relevant to check whether other features could allow us to achieve string stability without using absolute velocity in the controller, i.e.~with $h=0$.
\newline - The use of dynamic sensor mounts $M^{(r)}(s)$, $M^{(f)}(s)$ is another way of breaking the symmetry of purely relative position measurements $e_i = x_{i-1} - x_i$ when $h=0$. Tuning these dynamics in some beneficial way appears as a less invasive solution, than assuming with $h>0$ that a global reference is available for absolute velocity measurement. We thus see this as another potential \emph{opportunity} to break the impossibility observed in \cite{6,8}.
\newline - The use of local communication is a natural feature, considered e.g.~in \cite{7,17,18,19,20}, with $h>0$. We want to see how it fares with $h=0$. If communication was assumed to be perfect, it is tempting to impose $r_1 = x_0-x_1$, $r_i = r_{i-1} + e_i$ for $i>1$ such that in fact $r_i = x_0 - x_i$. It then becomes possible to control each $x_i$ individually with respect to the leader, on the basis of $u_i(s) = H(s) (x_0-x_i)$, as if we had a global reference for absolute position. In this case there is no distributed control problem anymore, and this is not the situation we want to consider. However, with imperfect communication as imposed by our model with Assumption 1, there is no way to obtain $r_i = x_0-x_i$ perfectly. Our model is thus just meant to crudely express communication uncertainty: banishing special deadbeat effects that would rely on perfect (pole) cancellation is a standard robustness assumption in linear control. Another typical limitation would be a finite bandwidth in $W$ associated with additive communication noise, but we do not even need this assumption for our results.


\section{Summary of results}\label{sec:3}

We now summarize more accurately the known results from the literature and new results of the present paper, see Table 1. Our proofs are given in appendix.

\begin{table*}
\caption{known results and \emph{clarifications provided by the present paper}; with $M^{(r)}=M^{(f)}=$identity unless specified.}
\begin{tabular}{|p{60mm}|p{60mm}|p{40mm}|}
\hline
Standard impossibility; \cite{6,8}; Lemma 1 &   $h=W=0$  &  Def.1 and Def.2 fail\\
\hline
Time-headway solution; \cite{15}; Proposition 1    &  $W=0$, $h>0$ with PD controller  &  Def.2 holds \\
\emph{Theorem 1}  & \emph{$W=0$, $h>0$ with $K(0)$ bounded} & \emph{Def.1 fails} \\
\emph{Theorem 2}  & \emph{$W=0$, $h>0$ with PID controller}  &  \emph{Def.1 and Def.2 hold} \\
\hline
CACC with time-headway; \cite{17,19} &  $W\neq0$, $h>0$  &  Def.2 holds with possibly better scaling, lower $h$\\
\emph{Theorem 3a}  & \emph{CACC-type $W\neq0$, $h>0$, $K(0)$ bounded}  &  \emph{Def.1 fails} \\
\emph{Theorem 3b}  & \emph{CACC type $W\neq0$, $h=0$}  &  \emph{Def.1 and Def.2 fail} \\
\emph{Theorem 4a}  & \emph{any $W\neq0$, $h=0$, $K(0)$ bounded}  &  \emph{Def.1 fails} \\
\emph{Theorem 4b}  & \emph{any $W\neq0$, scalar $v_i$, $h=0$}  &  \emph{Def.1 and Def.2 fail} \\
\hline
\emph{Theorem 5}  &  \emph{$W=0$, $h=0$, tuning $M^{(r)}$ and $M^{(f)}$} & \emph{Def.1 and Def.2 fail}\\
\hline
\end{tabular}
\end{table*}

Regarding known results, the first observations of string instability (\cite{8}) were made when each control input $u_i$ is reacting just to the relative distance $x_{i-1}-x_i$ with the vehicle in front, thus corresponding to \eqref{eq:2order},\eqref{eq:errmeas},\eqref{gc0} with $M^{(r)}=M^{(f)}=$identity, $h=W=0$. It was observed that $L_2$ string stability is impossible (Definition 2), and a fortiori $(L_2,l_2)$ string stability (Definition 1) as well, with any linear controller that avoids pole cancellation ($K(0) \neq 0$). Let us repeat the short argument, which has motivated alternatives. Let $u_i(s) = K(s) e_i(s)$. From \eqref{eq:2order}, the closed-loop equation for the $e_i$ and with a disturbance on the first subsystem only writes
$$e_i= T(s)^{i-1}\frac{1}{s^2+K(s)}d_0 \; ,$$
with $T(s)=\frac{K(s)}{s^2+K(s)} = \frac{R(s)}{1+R(s)}$ with $R(s)=K(s)/s^2$, taking the form of a complementary sensitivity function. To guarantee $L_2$ string stability, with $N$ unboundedly large, it is then necessary in particular that $\vert T(j\omega) \vert \leq 1$ at all frequencies $\omega$. This is impossible for a stable system, from the statement of Bode's Complementary Sensitivity integral, which we recall below.

\noindent \textbf{Lemma 1:} \emph{Assume that the loop transfer function $R(s)$ of a system has (at least) a double pole at $s=0$. If the associated feedback system is stable, then the complementary sensitivity function $T(s)=\frac{R(s)}{1+R(s)}$ must satisfy: 
$$\int_{0}^{\infty} ln\mid T(j\omega)\mid.d\omega/\omega^2 = \pi \sum_k\frac{1}{q^{(T)}_k} \ge 0 \, ,$$
where $\{ q^{(T)}_k \}$ are the zeros of $R(s)$ in the open right half plane. In particular, if $\vert T(j\omega) \vert < 1$ at some frequencies, then necessarily $\vert T(j\omega) \vert > 1$ at other frequencies.}

To circumvent this impossibility, \cite{15} have proposed a controller with $h>0$ (and still $W=0$ i.e.~no communication, $M^{(r)}=M^{(f)}=$identity i.e.~no special sensors). In that case the closed-loop dynamics writes
\begin{eqnarray}
\label{eq:controller2} e_i&=&\frac{K(s)}{s^2+(1+hs)K(s)}e_{i-1}\\ 
\nonumber & & +\frac{1}{s^2+(1+hs)K(s)}(d_{i-1}-(1+hs)d_i)
\end{eqnarray}
and we would have $T(s) = \frac{R(s)}{1+R(s)}$ with $R(s) = \frac{K(s)}{s^2 + h s K(s)}$, which does not have a double pole at $s=0$ and thus circumvents Lemma 1. For completeness and later comparison, we give the following result comparable to \cite{15}.

\noindent \textbf{Proposition 1:} \emph{The norm at $s=j\omega$ of transfer function $T(s)=\frac{K(s)}{s^2+(1+hs)K(s)}$ in \eqref{eq:controller2} is $< 1$ at all frequencies $\omega\neq0$, and its $H_\infty$ norm equals $T(0)=1$, if and only if one of the following equivalent conditions hold:
\newline (a) If one chooses $\bar{K}(s)=K(s)(1+hs)$ first and then derives $K(s)$ from $h$, then we should ensure that $h$ satisfies
\begin{equation}
\label{eq:headway1}
h>\sqrt{\underset{\omega}{\max}\frac{\left\vert\frac{\bar{R}(j\omega)}{1+\bar{R}(j\omega)}\right\vert^2-1}{\omega^2}}\
\end{equation}
in which $\bar{R}(s)=\bar{K}(s)/s^2$.
\newline (b) If one chooses $K(s)$ first, then the criterion becomes
\begin{equation}
\label{eq:headway2}
h>\underset{\omega}{\max} \sqrt{K_R(j\omega) \left(2-\omega^2 K_R(j\omega)\right)} + \omega K_J(j\omega)
\end{equation}
where $K_R(j\omega) = \frac{1}{2}(\frac{1}{K(j\omega)}+\frac{1}{K(j\omega)^*})$, $K_J(j\omega) = \frac{1}{2j}(\frac{1}{K(j\omega)}-\frac{1}{K(j\omega)^*})$, and the maximization runs over all $\omega$ for which the argument of the square root is positive.}

For particular controllers one can get easy criteria, e.g.~for a PD controller $K(s) = b s + a$, it is not hard to see that if $a>2b^2$ the right hand side in case (b) is decreasing with $\omega$, and the condition becomes $h>\sqrt{2/a}$.

A direct consequence of Prop.1 is that one can avoid amplifying a disturbance $d_0(s)$ along the vehicle chain, and satisfy Definition 2 of string stability. However, remarkably, a corresponding result about Definition 1, i.e.~the stronger version of string stability as considered in e.g.~\cite{2,9,11}, appears to be missing in the literature. The present paper provides the following results in this direction.

\noindent \textbf{Theorem 1:}  \emph{Consider the vehicle chain system \eqref{eq:2order},\eqref{eq:errmeas},\eqref{gc0} with $M^{(r)}=M^{(f)}=$identity and $W=0$. There exists no pair $\,(K(s),\;h)\,$, where $h\geq 0$ is any constant time-headway and $K(s)$ any stabilizing linear controller with $K(0)$ bounded, that would achieve $(L_2,l_2)$ norm string stability (Def.1).}

\noindent \textbf{Theorem 2:} \emph{Consider the vehicle chain system \eqref{eq:2order},\eqref{eq:errmeas},\eqref{gc0} with $M^{(r)}=M^{(f)}=$identity and $W=0$. A stabilizing PID controller $K(s)$ with headway $h$ satisfying \eqref{eq:headway2}, can ensure $(L_2,l_2)$ norm string stability (Def.1).}

In light of Thm.1, the lack of positive result about $(L_2,l_2)$ string stability in the time-headway literature can thus be attributed to their focus on \emph{bounded} stabilizing controllers $K(s)$, excluding e.g.~controllers with integral action. Theorem 2 further clarifies that a PID controller indeed does allow to achieve this stronger version of string stability.

The result of Thm.2 obviously covers as well the case where communication is allowed \emph{on top of time headway} $h>0$. This is the standard setting considered in papers like \cite{17,19}, where the communication is meant as a way to improve performance rather than just achieving string stability. Those papers also stay with the weaker notion of Def.2, and they impose a more precise communication structure, called Cooperative Adaptive Cruise Control (CACC), by assuming that the message sent by vehicle $i$ to its follower $i+1$ is a filtered version of the input command $u_i$:
\begin{equation}
\label{c2} v_{i}(s) = \tfrac{1}{B(s)} \big( K(s) e_{i}(s) + H(s) r_i(s) \big) \; ,
\end{equation}
or in other words \eqref{gc0} with $F = K/B$ and $G=H/B$.

By Thm.2, we have thus established that with $h>0$ and PID control, there exists a communication (namely the trivial one $F=H=G=0$) which does achieve the stronger string stability of Def.1 as well. A remaining question is then, \emph{how much can communication allow us to weaken some assumptions and still achieve string stability?} A negative answer is expressed by the following results.

\noindent \textbf{Theorem 3:} \emph{Consider the vehicle chain system \eqref{eq:2order},\eqref{eq:errmeas},\eqref{gc0} with $M^{(r)}=M^{(f)}=$identity, with the structure \eqref{c2} of CACC communication.
\newline \textbf{\emph{(a)}} There exists no $h>0$ and $W,B,H,K$ satisfying Assumption 1 with $K(0)$ bounded, allowing to satisfy $(L_2,l_2)$ string stability (Def.1).
\newline \textbf{\emph{(b)}} For $h=0$, there exists no $W,B,H,K$ satisfying Assumption 1 and allowing to satisfy $L_2$  (Def.2) and thus a fortiori $(L_2,l_2)$ string stability (Def.1).}

\noindent \textbf{Theorem 4:} \emph{Consider the vehicle chain system \eqref{eq:2order},\eqref{eq:errmeas},\eqref{gc0} with $M^{(r)}=M^{(f)}=$identity and $h=0$.
\newline \textbf{\emph{(a)}} For $K(0)$ bounded, there exist no $W,F,G,H,K$ satisfying Assumption 1, with each $v_i$ possibly an $n$-dimensional vector signal for some $n\geq1$, allowing to satisfy $(L_2,l_2)$ string stability (Def.1).
\newline \textbf{\emph{(b)}} For  $v_i$ scalar signals, there exist no $W,F,G,H,K$ satisfying Assumption 1 and allowing to satisfy $L_2$ (Def.2) and thus a fortiori $(L_2,l_2)$ string stability (Def.1).}

From Thm.3a, CACC-type communication together with time headway does not exempt us of requiring unbounded $K(0)$, i.e.~a result equivalent to Thm.1 still holds. Moreover, Theorems 3b and 4 indicate that communication does not allow us to avoid the necessity of time headway; most importantly, in a CACC-type setting or generalized scalar setting, even Def.2 cannot be satisfied. The assumptions a priori still leave a loophole for \emph{unbounded $K(0)$ and communicating vector signals $v_i$}, but we believe that this is just a technical issue and in fact we conjecture that no communication model satisfying Assumption 1 would allow to achieve string stability with \eqref{eq:2order},\eqref{eq:errmeas},\eqref{gc0} and $h=0$.

As a final attempt, we have checked whether we could replace the requirement of time headway $h>0$, by breaking the symmetry of relative position measurements via the dynamics $M^{(r)},M^{(f)}$. Unfortunately, here too the result is an impossibility.

\noindent \textbf{Theorem 5:} \emph{Consider the vehicle chain system \eqref{eq:2order},\eqref{eq:errmeas},\eqref{gc0} with $h=0$ and no communication ($W=G=F=H=0$). There exists no choice of stabilizing $K(s)$ and of $M^{(r)},M^{(f)}$ under the form \eqref{eq:sens2}, allowing to satisfy $L_2$ string stability (Def.2), and thus a fortiori $(L_2,l_2)$ string stability (Def.1).}


\section{Conclusion}\label{sec:conclusion} 

We have identified impossibilities to achieve string stability with linear controllers in several extended settings --- communication, sensor dynamics, time-headway with bounded DC gain --- and one possible solution, namely a PID controller with sufficient time headway. While the proofs do not involve complicated techniques, they do complete the picture about string stability in the strong sense ($(L_2,l_2)$ bounded reaction to simultaneous perturbations on all the vehicles) and weak sense (avoid amplification of disturbance on the leader). Regarding extended settings, the paper narrows down the options towards achieving string stability without using absolute velocity. Regarding bounded DC gain, the bad behavior close to zero frequency is a basic feature and versatile focus, which should allow to similarly complete the picture for a.o.~bidirectional controllers in future work. These results also highlight the necessity to carefully check in applications whether Def.1, Def.2 or possibly some other string stability notion is the right proxy for what really needs to be achieved. This point has indeed attracted surprisingly little discussion in the literature, although it does appear to change the possible conclusions.

We must end with a short outlook on open points. One main assumption in this line of literature, including our paper, is of course the use of linear models; nonlinearities and in particular quantization in digital controllers are a priori not covered by the tight impossibility results. We plan to address this, together with the remaining loophole of vector communication, in future work. Another point regarding communication, is the assumption of \emph{local} message transmission. If instead the vehicles were using a communication bus, the picture could be changed. Access to sending over the bus would have to be managed with event-driven decision logic, calling for a  comprehensive cyberphysical systems treatment of string stability. One would still have to investigate though, which breakthrough useful information such bus could transmit, when only local distances are measured.

\appendix


\section{Proofs}\label{sec:4}


\noindent \textbf{Proof of Proposition 1:}  For case (a), we reformulate $T(s)=\frac{1}{1+hs} \, \frac{\bar{K}(s)/s^2}{1+\bar{K}(s)/s^2}$. Then writing
$$|T(j\omega)|^2 = \frac{1}{1+\omega^2h^2}\, \left\vert\frac{R(j\omega)}{1+R(j\omega)}\right\vert^2 < 1 \;\;\; \text{for all } \omega\neq 0
$$
directly yields the expression, where the Bode integral (Lemma 1) ensures that the max inside the square root will be non-negative. For case (b), we just write $1/|T(j\omega)|^2 = |-\omega^2/K(j\omega) + (1 + h j \omega)|^2 > 1$ and we group real and imaginary parts to isolate $h$. \hfill $\square$\vspace{2mm}

\noindent \textbf{Proof of Theorem 1:}  We take in particular a disturbance input $d_0$ that affects the leading vehicle only. From \eqref{eq:controller2}, such disturbance leads to
\begin{eqnarray*}
e_i&=&T(s)^{i-1}\frac{1}{s^2+(1+hs)K(s)}d_0 \;\;\;, \;\; i=1,2,...,N,
\end{eqnarray*}
with $T(s)=\frac{K(s)}{s^2+(1+hs)K(s)}$ as in Proposition 1. Then
\begin{eqnarray}\label{new:thm1a}
\sum_{i=1}^{N} \vert e_i(s) \vert^2 &=& \sum_{i=0}^{N-1} \mid T(s)\mid^{2i} \cdot \frac{\vert d_0(s)\vert_2^2}{\vert s^2+(1+hs)K(s) \vert^2} \; .
\end{eqnarray}
Consider a disturbance concentrated at low frequencies, such that
$\int_{-\epsilon}^\epsilon \vert d_0(j\omega) \vert^2 d\omega \geq \frac{1}{2} \int_{-\infty}^{+\infty} \vert d_0(j\omega) \vert^2 d\omega \,$ for some $\epsilon \ll 1$.
Writing
\begin{equation}\label{new:thm1b}
\Vert e(.) \Vert_2^2 \geq  \int_{-\epsilon}^\epsilon \sum_{i=1}^{N} \vert e_i(j\omega) \vert^2\,   d\omega \;,
\end{equation}
we lower-bound the right hand side thanks to continuity at $\omega=0$. Indeed, in \eqref{new:thm1a}, for $K(0)$ finite, there exist $\delta,\alpha > 0$ such that $\frac{1}{\vert s^2+(1+hs)K(s) \vert^2} \vert_{s=j\omega} > \alpha$ for all $\omega \in (-\delta,\delta)$. Moreover, since $T(0)=1$, we can also make $\min_{\omega \in [0,\epsilon)} | T(j\omega) |$ arbitrarily close to $1$ by taking $\epsilon$ close enough to $0$. Therefore, we can make the right hand side of \eqref{new:thm1b} arbitrarily close to $\tfrac{N \alpha}{2} \; \Vert d_0(.) \Vert_2^2$ by concentrating $d_0$ on low enough frequencies $\epsilon<\delta$. Thus the factor relating $\Vert e(.) \Vert_2^2$ to $\Vert d_0(.) \Vert_2^2 = \Vert d(.) \Vert_2^2$ cannot be bounded independently of the disturbance signal $d(.)$ and of $N$. \hfill $\square$\vspace{2mm}

\noindent  \textbf{Proof of Theorem 2:} For this positive result we must prove that we can tune the gains such that the system is stable, and we can guarantee $\mid\mid e(.)\mid\mid_2<C_0 \, \mid\mid d(.)\mid\mid_2$ in which the constant $C_0$ is bounded independently of number of vehicles $N$.

For \emph{stability}: in the upcoming paragraphs, we will impose no particular tuning values to $h$ nor to the parameters of the PID controller. To satisfy stability, it is thus sufficient to find a PID controller and $h$ which make the system stable while fulfilling the conditions of Proposition 1. Considering the first criterion in Proposition 1, we will thus fixe some tuning of the polynomial $\bar{K}(s)$ which makes the system stable (just checking always the same denominator $s^3 + s (1+h\,s)\,K(s) = s^3 + s\,\bar{K}(s)$). Once $\bar{K}(s)$ has been selected, we would then choose $h$ according to the related criterion, while adapting the other parameters in order to maintain $\bar{K}(s)$ fixed as selected. For this to be possible, the only essential element is to prove that $h$ in the first criterion of Proposition 1 always remains bounded for a stable PID controller.

We thus consider $s^3 + s\,\bar{K}(s)$ to be any third-order polynomial with roots in the open left half plane. Then in the criterion,
$$\frac{\bar{R}}{1+\bar{R}} = \frac{s\,\bar{K}(s)}{s^3 + s\,\bar{K}(s)}$$
remains bounded for all $s=j\omega$ and we must only investigate the behavior for $\omega$ close to 0. From the inverse triangle inequality $|\frac{\bar{R}}{1+\bar{R}}|^2-1 \leq | (\frac{\bar{R}}{1+\bar{R}})^2 - 1|$ a sufficient criterion for Proposition 1(a) is
$$h >  \sqrt{\left\vert\frac{(j\omega)^2 \bar{K}^2(j\omega) - [(j\omega)^3 + (j\omega) \bar{K}(j\omega)]^2 }{\omega^2 [(j\omega)^3 + (j\omega) \bar{K}(j\omega)]^2}\right\vert} \; ,$$
which just comes down to 
$$h> \sqrt{\left\vert\frac{\omega^4 - 2 \omega^2 K(j\omega)}{(j\omega \bar{K}(j\omega)-j\omega^3)^2}\right\vert}\; .$$
For $\omega$ close to 0 and $K(j\omega)$ a PID controller, the dominating term is $h > \sqrt{\vert 2\omega\, k_I  / k_I^2 \vert}$, with $k_I$ the integral gain. This imposes a bounded constraint on $h$ and it is thus possible indeed to satisfy stability and the criterion of Proposition 1 simultaneously with a PID controller.\\

For \emph{string stability}, we write in matrix form:
$$ e(s) = \big(-L(s) \mathbf{A}  + L(s) \mathbf{B}(s) + P(s) \mathbf{C}(s)\big) \, d(s)$$
with $\; P(s)=\frac{s^2+hsK(s)}{(s^2+(1+hs)K(s))^2}\;$,
$\;L(s)=\frac{1}{s^2+(1+hs)K(s)}\;$ and the $N \times (N+1)$ matrices
\begin{eqnarray} 
    \mathbf{A} &=& \begin{bmatrix}\
    0 & 1 & 0 & \dots  & 0 \\
    0 & 0 & 1& \dots  & 0 \\
        0 & 0 & 0& \dots  & 0 \\
    \vdots & \vdots & \vdots & \vdots & \vdots \\
    \nonumber 0 &0 & 0 & \dots &    1\\
\end{bmatrix} , \;\;
    \mathbf{B}(s) = \begin{bmatrix}\
    1 & 0 & 0 & \dots  & 0 \\
    T(s) & 0 & 0& \dots  & 0 \\
        T(s)^2 & 0 & 0& \dots  & 0 \\
    \vdots & \vdots & \vdots & \vdots & \vdots \\
    \nonumber T(s)^{N-1} &0 & 0 & \dots &    0\\
\end{bmatrix}  \\[3mm]
    \mathbf{C}(s) &=&\begin{bmatrix}\
    0 & 0 & 0 & \dots  & 0 \\
    0 & 1 & 0& \dots  & 0 \\
        0 & T(s) & 1& \dots  & 0 \\
    \vdots & \vdots & \vdots & \vdots & \vdots \\
    \nonumber 0 &T(s)^{N-2} & T(s)^{N-3} & \dots &    0\\
\end{bmatrix} \;\; .
\end{eqnarray}    
By the triangle inequality,
\begin{eqnarray}\label{newthm2a}
\Vert e(s)\Vert_{2}&\le& \big( \vert L(s) \vert\, \Vert \mathbf{A}\Vert_2 + \vert L(s) \vert\, \Vert \mathbf{B}(s)\Vert_2 \\
\nonumber && \;\; + \vert P(s) \vert\, \Vert \mathbf{C}(s)\Vert_2 \; \big)  \; \Vert d(s) \Vert_2
\end{eqnarray}     
with the induced matrix norms
$\Vert \mathbf{D} \Vert_2 = \sqrt{\lambda_\text{max} (\mathbf{D}^* \mathbf{D})} \, $
where $^*$ is the complex conjugate transpose.
The proof now comes down to proving a bounded norm, independent of $N$ and $s=j\omega$, for each of the three terms in front of $\Vert d(s) \Vert_2$ in \eqref{newthm2a}.

For the first term, since $\mathbf{A}^*\mathbf{A} = \text{diag}(0,1,1,1,...,1)$, we immediately have $\vert L(s) \vert\, \Vert \mathbf{A}\Vert_2 = \vert L(s) \vert$, and the latter can be bounded independently of $s=j\omega$ for a stable system.

For the second term, we have $\mathbf{B}^* \mathbf{B} = \text{diag}\big( \frac{1-\vert T(s) \vert^{2N}}{1-\vert T(s)^2 \vert},\; 0,\;0,...,0\big)\, .$
Under the conditions of Prop.1, the numerator is lower than 1 and 
$\vert L(j\omega) \vert\, \Vert \mathbf{B}(j\omega)\Vert_2 \leq \sqrt{\frac{\vert L(j\omega)\vert^2}{1-\vert T(j\omega)\vert^2}} \; .$ The unbounded DC gain ensures that, when $T(j\omega)-1$ converges to 0 at $\omega=0$, so does $L(j\omega)$. Analysis close to $\omega = 0$ yields
$$\frac{\vert L(j\omega)\vert^2}{1-\vert T(j\omega)\vert^2} \simeq \frac{\omega^2}{k_I^2} \cdot \frac{1}{h^2\omega^2} = \frac{1}{k_I^2 h^2}$$
i.e.~the limit for $\omega \rightarrow 0$ of $\vert L(j\omega) \vert\, \Vert \mathbf{B}(j\omega)\Vert_2$ is bounded, independently of $N$. It is then easy to find a bound that is valid at all frequencies $\omega$, independently of $N$.

To bound the third term, we can use the Gershgorin disk theorem on the matrix $|P|^2\, \mathbf{C}^* \mathbf{C}$. For $|T(j\omega)| \leq 1$ we can bound finite sums of powers of $|T(j\omega)|$ by an infinite geometric series, and trivially check that the eigenvalues are bounded independently of $N$, for all $\omega$ outside a neighborhood of the origin $\omega=0$. The latter is indeed the only place where $|T(j\omega)| = 1$, and an expansion for $\omega \ll 1$ shows that in fact $\vert P(j\omega) \vert^2\, \Vert \mathbf{C}(j\omega)\Vert^2_2$ converges to zero for $\omega\rightarrow 0$. This concludes the proof. \hfill $\square$\vspace{2mm}

\noindent  \textbf{Proof of Theorem 3:} Defining $z_i = [e_i\; ;\; v_{i-1}-v_i ]$, the closed-loop dynamics with $d_0\neq 0$ only is described by:
\begin{eqnarray*}
&&z_{i+1} = \mathbf{T}(s)\, z_i \; , \;\;\; z_1 = 
\left[\begin{array}{l}
\frac{1}{s^2+K(1+hs)} \\ \frac{-K/B}{s^2+K(1+hs)}
\end{array} \right] \; d_{0} \phantom{KK}\\
&&  \text{with }\; \mathbf{T}(s)=\left[\begin{array}{ll}
\frac{K}{s^2+K(1+hs)} & \frac{HW}{s^2+K(1+hs)} \\
\frac{K}{B} \cdot\frac{s^2}{s^2+K(1+hs)} & \frac{HW}{B} \cdot \frac{s^2}{s^2+K(1+hs)}
\end{array} \right] \; .
\end{eqnarray*}
The key simplification implied by CACC is that $\mathbf{T}(s)$ is singular for all $s$, since the right column equals $HW/K$ times the left column. Thus the single nonzero eigenvalue of $T(s)$ equals its trace,
$\; \text{trace}(\mathbf{T}(s)) =  \frac{K+\frac{HW}{B}s^2}{K(1+hs) + s^2} \; .$

(a) For any $h>0$, $\text{trace}(\mathbf{T}(s))$ approaches $1$ when $s$ approaches 0. Thus like in the proof of Thm.1, a uniform bound over $\omega,N$ cannot be found for $\Vert e \Vert_2$, if the corresponding mode with $s=0$ has a nonzero component in $z_1$. The corresponding eigenvector at $s=0$ is $z_i \propto [1\; ; \;0]$, while the zero eigenvector is $z_i \propto [1 \; ; \; -K(0)/HW(0)]$. Thus $z_1$ will  have a nonzero component on the ``bad'' mode unless $z_1 \propto [1 \; ; \; -K/B(0)]\propto [1 \; ; \; -K(0)/HW(0)]$, i.e.~either $K(0)$ unbounded which we exclude by assumption, or $HW(0) = B(0)$ exactly. The latter is forbidden by the last requirement of Assumption 1.

(b) For $h=0$, we can rewrite $\text{trace}(\mathbf{T}(s)) = \frac{R}{1+R}$ with $\; R = \frac{K+\frac{HW}{B}s^2}{s^2(1-HW/B)} \;$. Since $\frac{HW}{B}$ is bounded, the denominator decays at least as $s^2$ for $s$ close to $0$ and we are in the conditions to apply Lemma 1 (Bode complementary sensitivity integral); this implies that there will be a range of frequencies $\omega$ where $\mathbf{T}(j\omega)$ has an eigenvalue with norm $|\text{trace}(\mathbf{T}(j\omega))|$ larger than $1$. As for case (a), with the last requirement of Assumption 1 the system will unavoidably have a component of $z_1(\omega)$ on this mode, which unavoidably makes the system string unstable in the sense of Def.2. \hfill $\square$

\noindent \textbf{Proof of Theorem 4:} Similarly to the case of CACC, and somewhat simplified thanks to $h=0$, by defining $z_i = [e_i \;;\; v_{i-1} ]$ we can reformulate the dynamics as:
\begin{eqnarray*}
&& z_{i+1} = \mathbf{T}(s)z_i  
\; + \; 
\left[\begin{array}{l}
\frac{1}{s^2+K} \\ 0
\end{array} \right] \; (d_{i}-d_{i+1})\\
&&\text{with }\;  \mathbf{T}(s)=\left[\begin{array}{ll}
 \frac{K-HWF}{s^2+K} & \frac{HW(I-GW)}{s^2+K} \\  F & GW 
\end{array} \right] \; .
   \end{eqnarray*} 
Here $I$ is the identity matrix, emphasizing that $v_i$ might be a vector and $F,G,H,W$ appropriate matrices.

For case (a) the proof follows the same lines as Thm.3(a), after checking that $\mathbf{T(0)}$ has an eigenvalue $1$. 

For case (b) we use a Routh-Hurwitz type criterion for discrete systems, see e.g.~\cite{21}. For a two-dimensional state matrix $A$, it states that the eigenvalues belong to the unit circle provided
\newline $\bullet$ $|\text{det}(A)| \leq 1  \quad \text{and}$
\newline $\bullet$ $\text{det}(A)^*\text{trace}(A) - \text{trace}(A)^*| \leq 1-|\text{det}(A)|^2 \; .$
\newline The determinant of $\mathbf{T}(s)$ imposes
$$|\text{det}| = \left\vert \frac{GWK-HWF}{s^2+K} \right\vert =: \frac{|A|}{|s^2+K|} \leq 1 \; ,$$
where we have defined $A = (GK-HF)W$. 
Next, we need
\begin{eqnarray*}
1 & \geq & \frac{|\text{trace} - \text{trace}^* \text{det}|}{1-|\text{det}|^2 } \\
&=& \left\vert 1 + \frac{s^2}{1-|\text{det}|^2} \,\left(\frac{GW-1}{s^2+K} - \frac{(GW-1)^*}{(s^2+K)^*} \frac{A}{s^2+K} \right)  \right\vert \; .
\end{eqnarray*}
Since $\frac{s^2}{1-|\text{det}|^2}$ is real negative for $s=j\omega$ and $|\text{det}|<1$, the above equation cannot be satisfied if $(GW-1)/(s^2+K)$ takes a real negative value for some $s=j\omega$. Indeed, for any $c_1,c_2$ real negative and $c_3$ complex but of norm smaller than one, we have that $1+ c_1 c_2 (1-c_3)$ lies outside the unit disk. Thus to conclude the proof, there remains to show that $(GW-1)/(s^2+K)$ will always take a real negative value for some $s=j\omega$.

Since $s^2+K$ has two more zeros than poles, and all zeros must satisfy stability, we have that the phase Bode plot of $1/(s^2+K)$ goes down at least by $180$ degrees, to end at $-180$ degrees for $\omega$ tending to infinity. In contrast, $GW-1$ has as many zeros as poles; all poles are stable, implying 90 degrees down in the phase Bode plot,  such that overall with $GW-1$ we either go down or stay, and again we end at $-180$ degrees for $\omega$ tending to infinity. Now assume as a first possibility, that $GW-1$ starts at another value than $-180$ degrees. In this case, it must go down nontrivially, i.e.~we must go sown by strictly more than 180 degrees to end up at $-360$ degrees: somewhere in between, there will be a 180 degree phase, proving impossibility. (Note indeed that we forbid any perfect cancellation with $GW=1$ at a target value of $\omega$.) So the only choice left is that $GW$ starts at -180 degrees. Then for $K(0)$ finite we would have a negative real phase at $s=0$, thus impossible. There remains the case with $K$ having a pole of order $m>0$ at $s=0$. In this case, $1/(s^2+K)$ has $m$ of its zeros at $s=0$, and $\frac{1}{s^m}\frac{1}{s^2+K}$ has a phase Bode plot going down by $(180+m 90)$ degrees overall. This means, $\frac{GW-1}{s^2+K}$ would start with a phase of $-180+m 90$ degrees at $s=0$, then go down by $180+m 90$ degrees to end up at $-360$ degrees for $\omega$ tending to infinity, with $m>0$. Again, this implies a phase of $-180$ degrees for some intermediate $\omega$. There are no possibilities left, so the proof is concluded.\hfill $\square$

\noindent \textbf{Proof of Theorem 5:} The error dynamics write
\begin{equation}\label{eq:snespr}
e_i = M^{(r)} \cdot \frac{1}{M^{(f)}} \cdot T'(s) \; e_{i-1} =: A(s) \; e_{i-1}
\end{equation}
where $T'(s)=\frac{K'(s)}{s^2+K'(s)}$ with $K'(s) = M^{(f)}(s) \cdot K(s)$.  With \eqref{eq:sens2}, $M^{(r)}$, $M^{(f)}$ and $T'$ all take the form of complementary sensitivity functions satisfying Lemma 1. We thus have
\begin{eqnarray*}
\int_0^{\infty} \text{ln}|A(j\omega)| . d\omega/\omega^2 & = &  \sum_k \frac{1}{q^{(M^{(r)})}_k} - \frac{1}{q^{(M^{(f)})}_k} + \frac{1}{q^{(T')}_k} \; ,
\end{eqnarray*}
with the $q$'s denoting the respective zeros of the loop transfer functions in the right half plane. Having a control effect requires $|A(j\omega)|<1$ at some frequencies, while having $|A(j\omega)|>1$ at any frequency would imply that the system is string unstable. Combining these two features requires that the right hand side be negative. The only way to obtain this is if $K^{(f)}(s)/s^2$ has zeros in the open right half plane, without having the same zeros in the other terms. However, the latter would mean that $M^{(f)}(s)$ has zeros in the right half plane, unmatched by the other transfer functions, and by \eqref{eq:snespr} this would imply that the vehicle chain has a pole in the right half plane i.e.~it is unstable.
\hfill $\square$

\end{document}